\documentclass[12pt,a4paper]{article} 
\usepackage[T1]{fontenc}
\usepackage{mathptmx,courier,pifont}
\usepackage[dvips]{graphicx}
\usepackage{psfrag}
\usepackage{amsmath,amssymb}        
\usepackage{latexsym}
\usepackage{url}
\def\RR{\hbox{I\kern-.2em\hbox{R}}}
\def\ds{\displaystyle}
\def\e{\epsilon}

\newcommand{\eqnsection}{
   \renewcommand{\theequation}{{\thesection.\arabic{equation}}}
   \makeatletter
   \csname @addtoreset\endcsname{equation}{section}
   \makeatother}
\eqnsection
\title{Blasius Problem and Falkner-Skan model: T{\"o}pfer's Algorithm and its Extension}
\author{Riccardo Fazio \\
Department of Mathematics and Computer Science \\
University of Messina \\
Viale F. Stagno D'Alcontres 31, 98166 Messina, Italy\\
{\small e-mail: rfazio@unime.it \ \ \ 
home-page: http://mat521.unime.it/fazio}}
\date{\small Submitted: October 18 and in revised form October 30 and December 15, 2012\\
Dedicated to the centenary of T{\"o}pfer's paper.}
\begin{document}               
\maketitle

\begin{abstract}
In this paper, we review the so-called T{\"o}pfer algorithm that allows us to find a non-iterative numerical solution of the Blasius problem, by solving a related initial value problem and applying a scaling transformation. 
Moreover, we remark that the applicability of this algorithm can be extended to any given problem, provided that the governing equation and the initial conditions are invariant under a scaling group of point transformations and that the asymptotic boundary condition is non-homogeneous. 
Then, we describe an iterative extension of T{\"o}pfer's algorithm that can be applied to a general class of problems. 
Finally, we solve the Falkner-Skan model, for values of the parameter where multiple solutions are admitted, and report original numerical results, in particular data related to the famous reverse flow solutions by Stewartson.
The numerical data obtained by the extended algorithm are in good agreement with those obtained in previous studies.
\end{abstract}

\bigskip
\bigskip



\section{Introduction}\label{Introduction}
At the beginning of the last century, Prandtl \cite{Prandtl:1904:UFK} put the foundations of
boundary-layer theory providing the basis for the unification of two sciences, which at that time seemed incompatible: namely, theoretical hydrodynamics and hydraulics.
The main application of boundary-layer theory are devoted to the calculation of the skin-friction drag acting on a body moving through a fluid, for example the drag of: an airplane wing, a turbine blade, or a complete ship (see Schlichting and Gersten \cite{Schlichting:2000:BLT}).
In this new century, due to the increasing number of applications of microelectronics devices, boundary-layer theory has found a renewal of interest within the study of gas and liquid flows at the micro-scale regime (see for instance Gad el Hak \cite{Gad-el-Hak:1999:FMM} or Martin and Boyd \cite{Martin:2001:BBL}).
Recently, Boyd \cite{Boyd:2008:BFC} has used the problem considered by Prandtl as an example where some good analysis allowed researchers of the past to solve problems, governed by partial differential equations, that might be otherwise impossible to face before the computer invention.
In this context, Prandtl, Blasius and T{\"o}pfer work is still modern now as it was more than one hundred years ago.

For an exhaustive derivation of the complete boundary layer equations governing general compressible flows, the interested reader is referred to Stewartson \cite{Stewartson:1964:TLB} or Schlichting and Gersten \cite{Schlichting:2000:BLT}.
Within boundary-layer theory, Blasius problem \cite{Blasius:1908:GFK} is given by:
\begin{eqnarray}
&{\displaystyle \frac{d^3 f}{d \eta^3}} + \frac{1}{2} f
{\displaystyle \frac{d^{2}f}{d\eta^2}} = 0 \nonumber \\[-1.5ex]
\label{eq:Blasius} \\[-1.5ex]
&f(0) = {\displaystyle \frac{df}{d\eta}}(0) = 0 \ , \qquad
{\displaystyle \frac{df}{d\eta}}(\eta) \rightarrow 1 \quad \mbox{as}
\quad \eta \rightarrow \infty \ , \nonumber
\end{eqnarray}
where $\eta$ and $f$ are suitable similarity variables.
This is a two point boundary value problem (BVP) obtained for the model describing the steady plane fluid flow past a thin flat plate.
 
Boyd has pointed out how this particular problem of boundary-layer theory has arose interest of prominent scientists, like H. Weyl, J. von Neumann, M. Van Dyke, etc.; see Table~1 in \cite{Boyd:1999:BFC}. 
The main reason for this interest is due to the hope that any approach developed for this epitome can be extended to more difficult hydrodynamics problems.

Blasius main interest was to compute, without worrying about existence or uniqueness of its BVP solution, the value of 
\begin{equation}\label{eq:skin-friction}
 \lambda = \frac{d^2f}{d\eta^2}(0) \ ,
\end{equation}
i. e., the skin-friction coefficient.
In order to compute this value, Blasius used a formal series solution around $\eta=0$ and an asymptotic expansions for large values of $\eta$, adjusting the constant $\lambda$ so as to connect both expansions in an intermediate region.
In this way, Blasius obtained the (wrong) bounds $ 0.3315 < \lambda < 0.33175$.

In 1912, T{\"o}pfer \cite{Topfer:1912:BAB} revised the work by Blasius and solved numerically the Blasius equation with initial conditions 
\begin{equation}\label{eq:TICs}
f(0)=\frac{df}{d\eta}(0)=0 \ , \qquad \frac{d^2f}{d\eta^2}(0)=1 \ . 
\end{equation}
Using hand computations with the classical fourth order Runge-Kutta method and a suitable scaling invariance, he arrived, without detailing his computations, at the value
$\lambda \approx 0.33206$, contradicting the bounds obtained by Blasius.  

Thereafter, the quest for a good approximation of $\lambda$ became a main concern.
This is seldom the case for the most important problems of applied mathematics:
the first study is usually devoted to find a method to solve a given problem, and once the problem is solved the attention is turned to know how accurate is the computed solution and whether there are different methods that can provide a solution with less effort or more accuracy. 
Using a power series, Bairstow \cite{Bairstow:1925:SF} reports $\lambda \approx 0.335$, whereas Goldstein \cite{Goldstein:1930:CSS} obtains $\lambda \approx 0.332$; moreover, using a finite difference method, Falkner \cite{Falkner:1936:MNS} computes $\lambda \approx 0.3325765$, and Howarth \cite{Horwarth:1938:SLB} yields $\lambda \approx 0.332057$.
Fazio \cite{Fazio:1992:BPF}, using a free boundary formulation of the Blasius problem, computes $\lambda \approx 0.332057336215$.
Boyd \cite{Boyd:1999:BFC} uses T{\"o}pfer's algorithm to obtain the accurate value $\lambda \approx 0.33205733621519630$.
By the Adomain's decomposition method, Abbasbandy \cite{Abbasbandy:2007:NSB} finds $\lambda \approx 0.333329$, whereas a variational iteration method with Pad\'e approximants allows Wazwaz \cite{Wazwaz:2007:VIM} to calculate the value $\lambda \approx 0.3732905625$.
Tajvidi et al. \cite{Tajvidi:2008:MRL} apply modified rational Legendre functions to get a value of  $\lambda \approx 0.33209$.

Also the Crocco formulation \cite{Crocco:1941:SLL} can be applied in order to compute the value of skin-friction coefficient, we can apply .
For instance, Vajravelu et al. \cite{Vajravelu:1991:SSS} use the Runge-Kutta method and a shooting technique to solve numerically the Crocco formulation and obtain the value $\lambda \approx 0.3322$, whereas Callegari and Friedman \cite{Callegari:1968:ASN} reformulate the Blasius problem in terms of the Crocco variables, show that this problem has an analytical solution, and compute the following bounds: $0.332055 < \lambda < 0.33207$.

Within scaling invariance theory, the transformation of a BVP into an initial value problem (IVP) due to T{\"o}pfer is a consequence of the invariance of Blasius equation and the two boundary conditions at $\eta = 0$ (at the plate) with respect to the scaling group of transformations  
\begin{equation}\label{eq:gscaleinv}
f^* = \lambda^{-\alpha} f \ , \qquad \eta^* = \lambda^{\alpha} \eta  \ ,   
\end{equation}
where $ \alpha $ is a non-zero parameter.
Using this transformation, a simple existence and uniqueness theorem was given by J. Serrin, as reported by Meyer \cite[pp. 104-105]{Meyer:1971:IMF}. 
This scaling invariance is essential also to the error analysis of the truncated boundary formulation for Blasius problem due to Rubel \cite{Rubel:1955:EET}.
A more complex proof of the existence and uniqueness of the Blasius problem solution was given before by Weyl \cite{Weyl:1942:DES}, who also proved that the solution has a positive second order derivative that is monotone decreasing on $ [0, \infty) $ and approaches to zero as $ \eta $ goes to infinity.

Our main goal here is to show how to solve numerically Blasius problem and similar problems in boundary layer theory, by initial value methods defined within scaling invariance theory.
These methods are referred to as numerical transformation methods (TMs).
As pointed out by Na \cite[Chapters 7-9]{Na:1979:CME}, usually a given, even simple, extension of the Blasius problem cannot be solved by T\"opfer algorithm.
Therefore, in order to extend the applicability of this non-iteratve transformation method (ITM) an iterative
version has been developed in \cite{Fazio:1994:FSE,Fazio:1996:NAN}.

In particular, we apply the ITM to the Falkner-Skan equation with relevant boundary conditions
\begin{eqnarray}\label{eq:abf1}
& {\displaystyle \frac{d^{3}f}{d\eta^3}} + f 
{\displaystyle \frac{d^{2}f}{d\eta^2}} + \beta \left[ 1 - \left({\displaystyle
\frac{df}{d\eta}}\right)^2 \right] = 0 \nonumber \\[-1.2ex]
& \\[-1.2ex]
& f(0) = {\displaystyle \frac{df}{d\eta}}(0) = 0 \ , \qquad
{\displaystyle \frac{df}{d\eta}}(\eta) \rightarrow 1 \quad \mbox{as}
\quad \eta \rightarrow \infty \ , \nonumber
\end{eqnarray}
where $ f $ and $ \eta $ are appropriate similarity variables
and $ \beta $ is a parameter. 
This set is called the Falkner-Skan model, after the names of two English mathematicians
who first studied it \cite{Falkner:1931:SAS}.
As pointed out by Na \cite[pp. 146-147]{Na:1979:CME}, if $ \beta \ne 0$, the BVP (\ref{eq:abf1}) cannot be solved by a non-ITM.
Indeed, the governing differential equation in (\ref{eq:abf1}) is not invariant with respect to any scaling group of point transformations.

The existence and uniqueness question for the problem (\ref{eq:abf1}) is really a complex matter.
Assuming that $ \beta > 0 $ and under the restriction $ 0 < \frac{df}{d\eta} < 1 $, known as normal flow condition, Hartree \cite{Hartree:1937:EOF} and Stewartson \cite{Stewartson:1954:FSF} proved that the problem (\ref{eq:abf1}) has a unique solution, whose first derivative tends to one exponentially.
Coppel \cite{Coppel:1960:DEB} and Craven and Peletier \cite{Craven:1972:USF} pointed out that the above restriction on the first derivative can be omitted when $ 0 \le \beta \le 1$.
Weyl proved, in \cite{Weyl:1942:DES}, that for each value of the parameter $\beta$ there exists a physical solution with positive  monotone decreasing, in $[0, \infty)$, second derivative that approaches zero as the independent variables goes to infinity. 
In the case $ \beta > 1 $, the Falkner-Skan model loses the uniqueness and a hierarchy of solutions with reverse flow exists.
In fact, for $ \beta > 1$ Craven and Peletier \cite{Craven:1972:RFS} computed solutions for which $ \frac{df}{d\eta} < 0 $ for some value of $ \eta $.
In each of these solutions the velocity approaches its limit exponentially in $\eta$. 
As mentioned before, the term normal flow indicates that
the flow velocity has a unique direction, and instead, reverse flow means that the velocity
is both positive and negative in the integration interval. 

The considered problem has also multiple solutions for $\beta_{\min} < \beta < 0$, as reported by  Veldman and van de Vooren \cite{Veldman:1980:GFS}, with the minimum value of $\beta$ given by 
\begin{equation}\label{eq:bmin}
\beta_{\min} = -0.1988\dots \ .
\end{equation} 
In this range there exist two physical solutions, one for normal flow and one
for reverse flow.  
For $\beta = \beta_{\min}$ only one solution exists.
Finally, for $\beta < \beta_{\min}$ the problem has no solution at all. 
Our interest here is to apply the ITM to the range of $\beta$ where multiple solutions are admitted, in particular, our interest is to get numerically the famous solutions of Stewartson \cite{Stewartson:1954:FSF,Stewartson:1964:TLB}.
The obtained results are original and, as we shall see in the following sections, are in agreement with those available in literature.

The first computational treatment of the Falkner-Skan model is due to Hartree \cite{Hartree:1937:EOF}.
Cebeci and Keller \cite{Cebeci:1971:SPS} apply shooting and parallel shooting methods requiring asymptotic boundary condition to be imposed at a changing unknown boundary in the computation process. 
As a result, they report convergence difficulties, which can be avoided by moving towards more complicated methods. 
Moreover, to guarantee reasonable accuracy, they are forced to use a {\it small enough} step-size and extensive computation for the solution of the IVPs.
Na \cite[pp. 280-286]{Na:1979:CME} describes the application of invariant imbedding.
A modified shooting method \cite{Asaithambi:1997:NMS} and finite-difference methods \cite{Asaithambi:1998:FDM,Asaithambi:2004:SOF} for this problem are presented by Asaithambi.
Kuo \cite{Kuo:2003:ADT} uses a differential transformation method, which obtains a series solution of the Falkner-Skan equation.
Sher and Yakhot \cite{Sher:2001:NAS} define a new approach to solve this problem by shooting from infinity, using some simple analysis of the asymptotic behaviour of the solution at infinity.
Asaithambi \cite{Asaithambi:2005:SFE} proposes a faster shooting method by using recursive evaluation of Taylor coefficients.
Zhang and Chen \cite{Zhang:2009:IMS} investigate a modification of the shooting method, where the computation of the Jacobian matrix is obtained by solving two IVPs. 
A Galerkin-Laguerre spectral method is defined by Auteri and Quartapelle \cite{Auteri:2012:GLS}.

\section{T{\"o}pfer transformation}\label{Theory}
In order to clarify T{\"o}pfer \cite{Topfer:1912:BAB} derivation of a transformation of variables that reduces the BVP (\ref{eq:Blasius}) into an IVP we will consider the derivation of the series expansion solution.
Of course, some of the coefficients of the series can be evaluated by imposing the boundary conditions at $ \eta = 0 $.
Moreover, for the missing initial condition, we set
\begin{equation*}
     \frac{d^2f}{d\eta^2}(0) = \lambda \ , 
\end{equation*} 
where $ \lambda $ is a non-zero constant.
So that, we look for a series solution defined as
\begin{equation*}
    f(\eta) = \frac{\lambda}{2} \eta^2 + \sum_{n=3}^\infty C_n \eta^n
\end{equation*}
where the coefficients $ \lambda $ and $ C_n $, for $ n = 3, 4, \dots $, are constants to be determined.
In fact, the boundary values at the plate surface, at $\eta = 0$,
require that $ C_0 = C_1 = 0 $, and we also have $ C_2 = \lambda/2 $ by the definition of $ \lambda $.
Now, we substitute this series expansion into the governing differential equation, whereupon we find 
\begin{multline*} 
\sum_{n=3}^\infty n (n-1)(n-2)C_n \eta^{n-3} +\\
 + \frac{1}{2} 
\left(\frac{\lambda}{2} \eta^2 + \sum_{n=3}^\infty C_n \eta^n\right)
\left[\lambda + \sum_{n=3}^\infty n (n-1) C_n \eta^{n-2} \right] = 0 
\end{multline*}
or in expanded form 
\begin{multline*}
 \left[3 \cdot 2 \cdot C_3 \right] 
+ \left[4 \cdot 3 \cdot 2 \cdot C_4 \right] \eta 
+ \left[5 \cdot 4 \cdot 3 \cdot C_5 + \frac{1}{2} \cdot 2 \cdot \frac{\lambda}{2} \frac{\lambda}{2}\right] \eta^2 + \\
 + \left[6 \cdot 5 \cdot 4 \cdot C_6 + \frac{1}{2} \cdot 2 \cdot \frac{\lambda}{2} C_3 
 + \frac{1}{2} \cdot \frac{\lambda}{2} \cdot 3 \cdot 2 \cdot C_3\right] \eta^3 
 + \cdots = 0 \ .
\end{multline*}
According to a standard approach we have to require that all coefficients of the powers of $ \eta $ to be zero.
It is a simple matter to compute the coefficients of the series expansion in terms of $ \lambda $:
\begin{align*}
C_3 &= C_4 = 0 \ , \qquad C_5 = -\frac{\lambda^2}{2\cdot 5!} \nonumber \\
C_6 &= C_7 = 0 \ , \qquad C_8 = 11\frac{\lambda^3}{2^2\cdot 8!} \nonumber \\
C_9 &= C_{10} = 0 \ , \qquad C_{11} = -375\frac{\lambda^4}{2^3\cdot 11!} \nonumber \\
&\mbox{and so on ...} \nonumber 
\end{align*}
The solution can be written as
\begin{equation*}
f = \frac{\lambda \eta^2}{2} - \frac{\lambda^2 \eta^5}{2\cdot 5!} 
+ \frac{11 \cdot \lambda^3 \eta^8}{2^2\cdot 8!} 
- \frac{375 \cdot \lambda^4 \eta^{11}}{2^3\cdot 11!} + \cdots 
\end{equation*}
where the only unknown constant is $ \lambda $.
In principle, $ \lambda $ can be determined by imposing the boundary condition at the second point, but in this case this cannot be done because the left boundary condition is given at infinity.
However, by modifying the powers of $ \lambda $ we can rewrite the series expansion as
\begin{multline*} 
\lambda^{-1/3} f = \frac{\left(\lambda^{1/3} \eta\right)^2}{2} - \frac{\left(\lambda^{1/3} \eta\right)^5}{2\cdot 5!} 
+ \frac{11 \cdot \left(\lambda^{1/3} \eta\right)^8}{2^2\cdot 8!} 
- \frac{375 \cdot \left(\lambda^{1/3} \eta\right)^{11}}{2^3\cdot 11!} + \cdots 
\end{multline*} 
which suggests a transformation of the form 
\begin{equation}\label{eq:scalinv:Blasius}
f^* = \lambda^{-1/3} f \ , \qquad \eta^* = \lambda^{1/3} \eta  \ .   
\end{equation}
In the new variables the series expansion becomes 
\begin{equation*}
f^* = \frac{\eta^{*2}}{2} - \frac{\eta^{*5}}{2\cdot 5!} 
+ \frac{11 \cdot \eta^{*8}}{2^2\cdot 8!} 
- \frac{375 \cdot \eta^{*11}}{2^3\cdot 11!} + \cdots 
\end{equation*}
which does not depend on $ \lambda $.
We notice that the governing differential equation and the initial conditions at the free surface, at $\eta = 0$, are left invariant by the new variables defined above.
Moreover, the first and second order derivatives transform in the following way
\begin{equation*}
\frac{d f^*}{d \eta^{*}} = \lambda^{-2/3} \frac{d f}{d \eta} \ ,   
\qquad
\frac{d^2 f^*}{d \eta^{*2}} = \lambda^{-1} \frac{d^2 f}{d \eta^{2}} \ .  
\end{equation*}
As a consequence of the definition of $ \lambda $ we have
\begin{equation*}
\frac{d^2 f^*}{d \eta^{*2}} (0) = 1 \ ,   
\end{equation*} 
and this explain why in these variables the series expansion does not depend on $ \lambda $.
Furthermore, the value of $ \lambda $ can be found provided that we have an approximation for $ \frac{d f^*}{d \eta^{*}}(\infty) $.
In fact, by the above relation we get
\begin{equation}\label{eq:lambda:Blasius}
\lambda = \left[ \frac{d f^*}{d \eta^{*}}(\infty) \right]^{-3/2} \ .   
\end{equation} 

From a numerical viewpoint, BVPs must be solved within the computational domain simultaneously (a \lq \lq stationary\rq \rq \ problem), whereas IVPs can be solved by a stepwise procedure (an \lq \lq evolution\rq \rq \ problem).
Somehow, numerically, IVPs are simpler than BVPs.  

\subsection{T{\"o}pfer algorithm}
Let us list the steps necessary to solve the Blasius problem by the T{\"o}pfer algorithm:
\begin{enumerate}
    \item we solve the IVP
\begin{eqnarray}
&{\displaystyle \frac{d^3 f^*}{d \eta^{*3}}} + \frac{1}{2} f^*
{\displaystyle \frac{d^{2}f^*}{d\eta^{*2}}} = 0 \nonumber \\[-1.5ex]
\label{eq:Blasius2} \\[-1.5ex]
&f^*(0) = {\displaystyle \frac{df^*}{d\eta^*}}(0) = 0, \qquad
{\displaystyle \frac{d^2f^*}{d\eta^{*2}}}(0) = 1  \nonumber
\end{eqnarray}
and, in particular, get an approximation for
$ \frac{d f^*}{d \eta^{*}}(\infty) $ in order to
compute $ \lambda $ by equation (\ref{eq:lambda:Blasius});
\item
we obtain $ f(\eta) $, ${\displaystyle \frac{df}{d\eta}}(\eta)$, and ${\displaystyle \frac{d^2f}{d\eta^2}}(\eta)$ by the inverse transformation of (\ref{eq:scalinv:Blasius}).
\end{enumerate}
In this way we define a non-ITM.

\begin{figure}[p]
	\centering
\psfrag{e}[][]{$ \eta^*, \eta $} 
\psfrag{f2star}[][]{$ \frac{df^*}{d\eta^*} $} 
\psfrag{f2}[][]{$ \frac{df}{d\eta} $} 
\psfrag{f3star}[][]{$ \frac{d^2f^*}{d\eta^{*2}} $} 
\psfrag{f3}[][]{$ \frac{d^2f}{d\eta^2} $} 
\includegraphics[width=.7\textwidth]{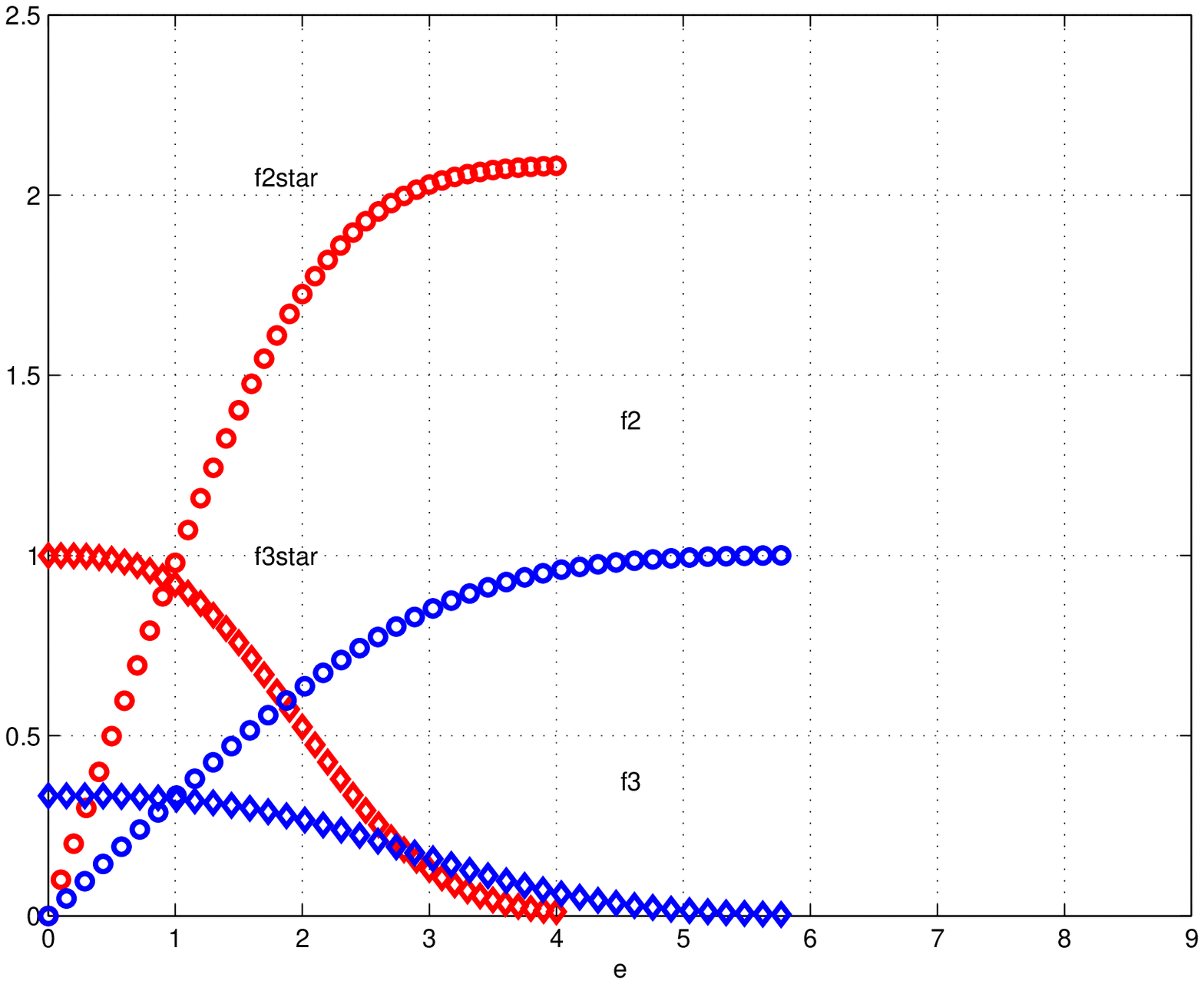} \\
\includegraphics[width=.7\textwidth]{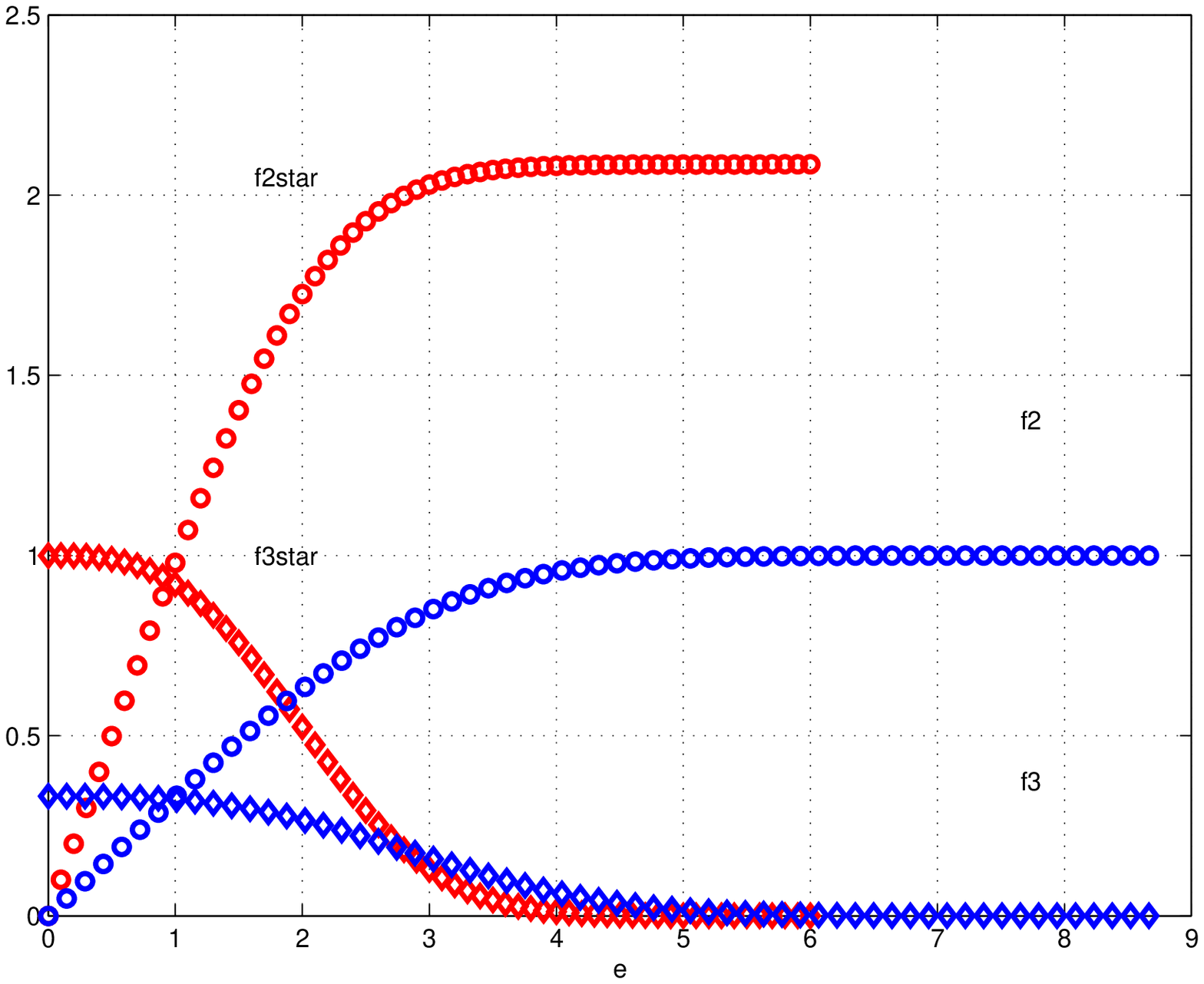}
\caption{Blasius solution solved with T\"opfer's algorithm. Top frame: for $\eta^* \in [0,4]$ we get $\lambda_1 = 0.3329124105$. Bottom frame: for $\eta^* \in [0,6]$ we find $\lambda_2 = 0.3320575595$.} 
	\label{fig:Blasius}
\end{figure}

Indeed, T{\"o}pfer solved the IVP for the Blasius equation once.
At large but finite $ \eta_j^* $, ordered so that $ \eta_j^* < \eta_{j+1}^* $, we can compute by equation (\ref{eq:lambda:Blasius}) the corresponding $ \lambda_j $.  
If two subsequent values of $ \lambda_j $ agree within a specified accuracy, then $ \lambda $ is approximately equal to the common value of the $ \lambda_j $, otherwise, we can march to a larger value of $ \eta $ and try again.
Using the classical fourth order Runge-Kutta method (see Butcher \cite[p. 166]{Butcher}) and a grid step $ \Delta \eta^* = 0.1$
T{\"o}pfer was able to determine $ \lambda $ with an error less than $ 10^{-5} $. 
He used the two truncated boundaries $\eta_1^* = 4$ and $\eta_2^* = 6$. 
We reproduce T{\"o}pfer computations in figure~\ref{fig:Blasius}.
In fact, in figure~\ref{fig:Blasius} we plot the numerical solutions obtained by T{\"o}pfer's algorithm defined above.
We notice that the top and the bottom frames of this figure show the solutions of the IVP (\ref{eq:Blasius2}) and of the BVP (\ref{eq:Blasius}).

\section{Extension of T{\"o}pfer algorithm}
The applicability of a non-ITM to the Blasius problem is a consequence of its partial invariance with respect to the transformation (\ref{eq:scalinv:Blasius}); the asymptotic boundary condition is not invariant.
The non-iterative algorithm can be extended, to a given problem in boundary layer theory, provided that the governing equation and the initial conditions are invariant under a scaling group of point transformations and the asymptotic boundary condition is non-homogeneous. 
Several problems in boundary-layer theory lack this kind of invariance and cannot be solved by non-ITMs
\cite[Chapters 7-9]{Na:1979:CME}. 
To overcome this drawback, we can introduce an iterative extension of the algorithm.
The main idea for the new algorithm is to modify the original problem, by introducing a numerical parameter $ h $, and to require the invariance of the modified problem with respect to an extended scaling group involving $ h $; see \cite{Fazio:1994:FSE,Fazio:1996:NAN} for details. 
\subsection{The iterative transformation method}
In order to define the ITM, let us consider the class of BVPs defined by
\begin{eqnarray}
&{\displaystyle \frac{d^3 f}{d \eta^3}} = \phi\left(\eta, f,
{\displaystyle \frac{df}{d\eta}, \frac{d^{2}f}{d\eta^2}}\right) \nonumber \\[-1.5ex]
\label{eq:class} \\[-1.5ex]
&f(0) = a \qquad  {\displaystyle \frac{df}{d\eta}}(0) = b \ , \qquad
{\displaystyle \frac{df}{d\eta}}(\eta) \rightarrow c \quad \mbox{as}
\quad \eta \rightarrow \infty \ , \nonumber
\end{eqnarray}
where $a$, $b$ and $c \ne 0$ are given constants.
We modify the class of problems (\ref{eq:class}) by introducing a numerical parameter $ h $ as follows
\begin{eqnarray}
&{\displaystyle \frac{d^3 f}{d \eta^3}} = h^{(1-3 \delta)/\sigma}\phi\left(h^{-(\delta/\sigma)} \eta, h^{-1/\sigma} f,
{\displaystyle h^{(\delta-1)/\sigma} \frac{df}{d\eta}, h^{(2 \delta-1)/\sigma} \frac{d^{2}f}{d\eta^2}}\right) \nonumber \\[-1.5ex]
\label{eq:class:mod} \\[-1.5ex]
&f(0) = h^{1/\sigma} a \qquad  {\displaystyle \frac{df}{d\eta}}(0) = h^{(1-\delta)/\sigma} b \ , \qquad
{\displaystyle \frac{df}{d\eta}}(\eta) \rightarrow c \quad \mbox{as}
\quad \eta \rightarrow \infty \ . \nonumber
\end{eqnarray}
It is worth noticing that the extended problem (\ref{eq:class:mod}) reduces to the original problem  (\ref{eq:class}) for $h=1$.
Moreover, the extended problem (\ref{eq:class:mod}) is partially invariant, the asymptotic boundary condition is not invariant, with respect to the extended scaling group of transformations
\begin{equation}\label{eq:scaling}
f^* = \lambda f \ , \qquad \eta^* = \lambda^{\delta} \eta \ , \qquad 
h^* = \lambda^{\sigma} h \ ,   
\end{equation}
with $\delta \ne 1$ and $\sigma \ne 0$.
Therefore, to find a solution of the given BVP means to find a zero of the so-called {\it transformation function} 
\begin{equation}\label{eq:AA2.9}
\Gamma (h^{*}) = \lambda^{-\sigma} h^* - 1 \ , 
\end{equation}  
where the group parameter $ \lambda $ is defined by the formula
\begin{equation}\label{eq:lambda2}
\lambda = \left[\displaystyle \frac{df^*}{d\eta^*}(\eta_\infty^*)/c\right]^{1/(1-\delta)} \ ,
\end{equation}  
and to this end we can use a root-finder method.
Let us notice that $\lambda$ and the transformation function are defined implicitly by the solution of the IVP
\begin{eqnarray}
&{\displaystyle \frac{d^3 f^*}{d \eta^{*3}}} = h_j^{*(1-3 \delta)/\sigma}\phi\left(h_j^{*-(\delta/\sigma)} \eta, h_j^{*-1/\sigma} f,
{\displaystyle h_j^{*(\delta-1)/\sigma} \frac{df^*}{d\eta^*}, h_j^{*(2 \delta-1)/\sigma} \frac{d^{2}f^*}{d\eta^{*2*}}}\right) \nonumber \\[-.5ex]
\label{eq:ITM} \\[-1ex]
&f^*(0) = h_j^{*1/\sigma} a \qquad  {\displaystyle \frac{df^*}{d\eta^*}}(0) = h_j^{*(1-\delta)/\sigma} b \ , \qquad
{\displaystyle \frac{d^2f^*}{d\eta^{*2}}}(0) = d \ , \nonumber
\end{eqnarray}
where $d$ is a parameter fixed by the user.
In particular, we are interested to compute $\frac{df^*}{d\eta^*}(\eta_\infty^*)$, an approximation of the asymptotic value $\frac{df^*}{d\eta^*}(\infty)$, which is used in the definition of $\lambda$. 

We set the values of $d$ and $\sigma$ and follow the steps:
\begin{enumerate}
\item{} we apply a root-finder method to compute a sequence $ h^{*}_{j}, $ for $ j = 0, 1, 2, \dots, $ .
Two sequences $\lambda_j$ and $ \Gamma (h^{*}_{j}) $ for 
$ j = 0, 1, 2, \dots, \ $ are defined by equation (\ref{eq:lambda2}) and (\ref{eq:AA2.9}), respectively.
\item{} a suitable convergence criterion should be used to verify whether $ \Gamma (h^{*}_{j}) \rightarrow  0 $ as $ j \rightarrow \infty $.
In this way we find the correct value of $h^*$ that transforms into $ h = 1 $.
If this is the case, then $\lambda_j$ converges to the correct value of $\lambda$ in the same limit.
\item{} once the correct value of $\lambda$ has been found, the solution of the original problem can be obtained by rescaling.
In particular, we have that
\[
{\displaystyle \frac{d^2f}{{d\eta}^2}}(0) = \lambda^{2 \delta - 1}
{\displaystyle \frac{d^2f^*}{{d\eta^*}^2}}(0) \ .
\]
\end{enumerate}
In this way we define an ITM.

In the next sub-section we apply the above iterative extension of T{\"o}pfer algorithm to the Falkner-Skan model. 

\subsection{The Falkner-Skan model}
In order to apply an ITM to (\ref{eq:abf1}) we have to embed it to a modified model and require the invariance of this last model with respect to an extended scaling group of transformations.
This can be done in several ways that are all equivalent.
In fact, the modified model can be written as 
\begin{eqnarray}\label{eq:abf2}
& {\displaystyle \frac{d^{3}f}{d\eta^3}} + f 
{\displaystyle \frac{d^{2}f}{d\eta^2}} + \beta \left[ h^{4/\sigma} - \left({\displaystyle
\frac{df}{d\eta}}\right)^2 \right] = 0 \ , \nonumber \\[-1.2ex]
& \\[-1.2ex]
& f(0) = {\displaystyle \frac{df}{d\eta}}(0) = 0 \ , \qquad
{\displaystyle \frac{df}{d\eta}}(\eta) \rightarrow 1 \quad \mbox{as}
\quad \eta \rightarrow \infty \ , \nonumber
\end{eqnarray}
and the related extended scaling group is given by
\begin{equation}\label{eq:scalinv}
f^* = \lambda f \ , \qquad \eta^* = \lambda^{-1} \eta \ , \qquad 
h^* = \lambda^{\sigma} h \ ,   
\end{equation}
where $\sigma$ is a parameter.
In the following we set $\sigma = 4$; for the choice $\sigma=8$ see \cite{Fazio:1994:FSE}.
In \cite{Fazio:1994:FSE}, a free boundary formulation of the Falkner-Skan model was considered and
numerical results were computed for the Homann flow ($ \beta = 1/2 $) as well as for the Hiemenz flow ($ \beta = 1 $).

From a numerical point of view the request to evaluate 
$ \frac{d f}{d \eta} (\infty) $ cannot be fulfilled.
Several strategies have been proposed in order to provide an approximation of this value.
The simplest and widely used one is to introduce, instead of infinity, a suitable truncated boundary.
A recent successful way to deal with such a issue is to reformulate the considered problem as a free BVP \cite{Fazio:1992:BPF,Fazio:1994:FSE,Fazio:1996:NAN}; for a survey on this topic see \cite{Fazio:2002:SFB}. 
Recently, Zhang and Chen \cite{Zhang:2009:IMS} have used a free boundary formulation to compute the normal flow solutions of the Falkner-Skan model in the full range $ \beta_{\min} < \beta \le 40$.
They applied a modified Newton's method to compute both the initial velocity and the free boundary.  
For the sake of simplicity, we do not use the free boundary approach but,  following T\"opfer,  we use some preliminary computational tests to find a suitable value for the truncated boundary.

\begin{figure}[p]
	\centering
\psfrag{h*}[][]{$h^*$} 
\psfrag{G}[][]{$\Gamma(h^*)$} 

\includegraphics[width=.475\textwidth]{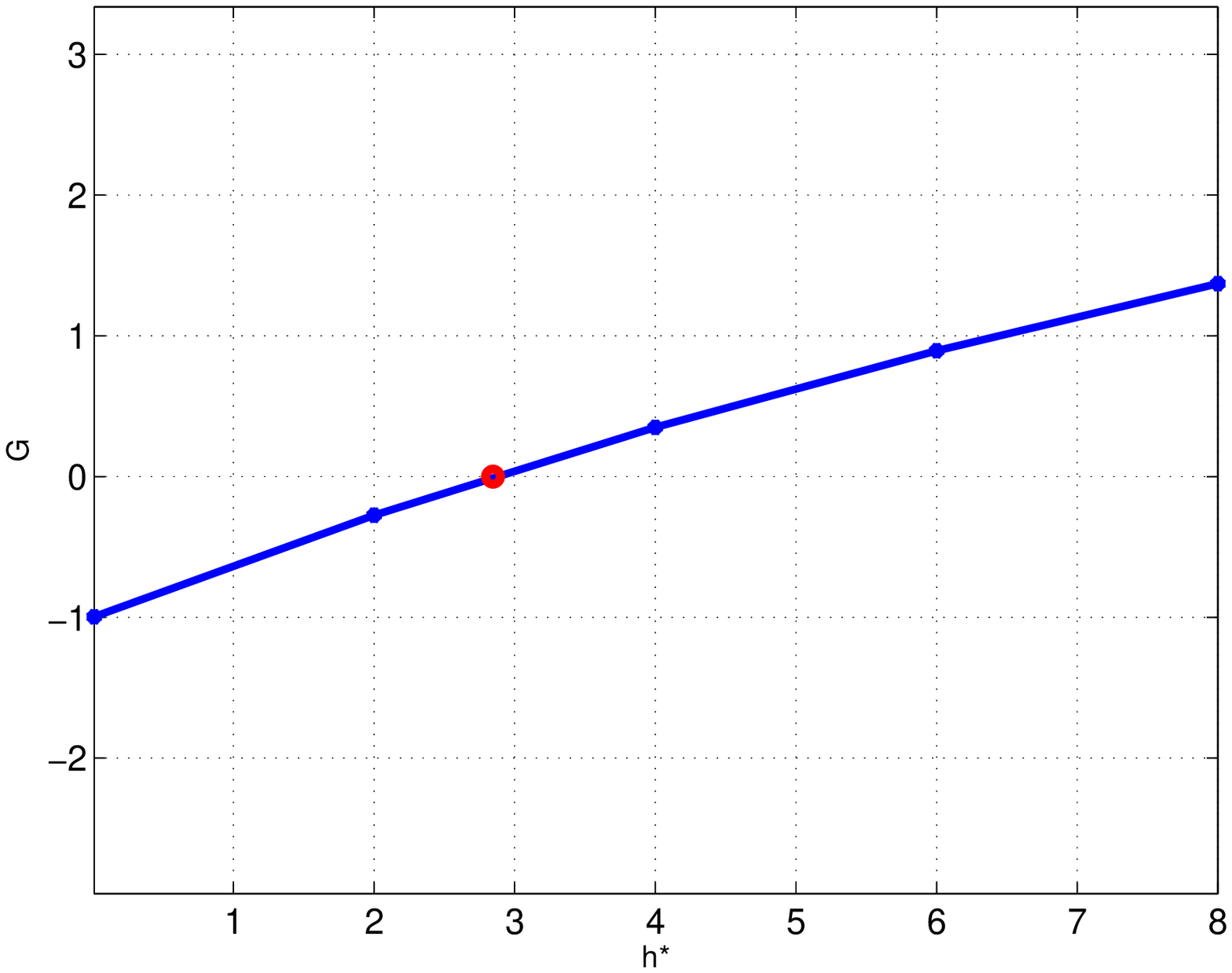}
\hfill
\includegraphics[width=.475\textwidth]{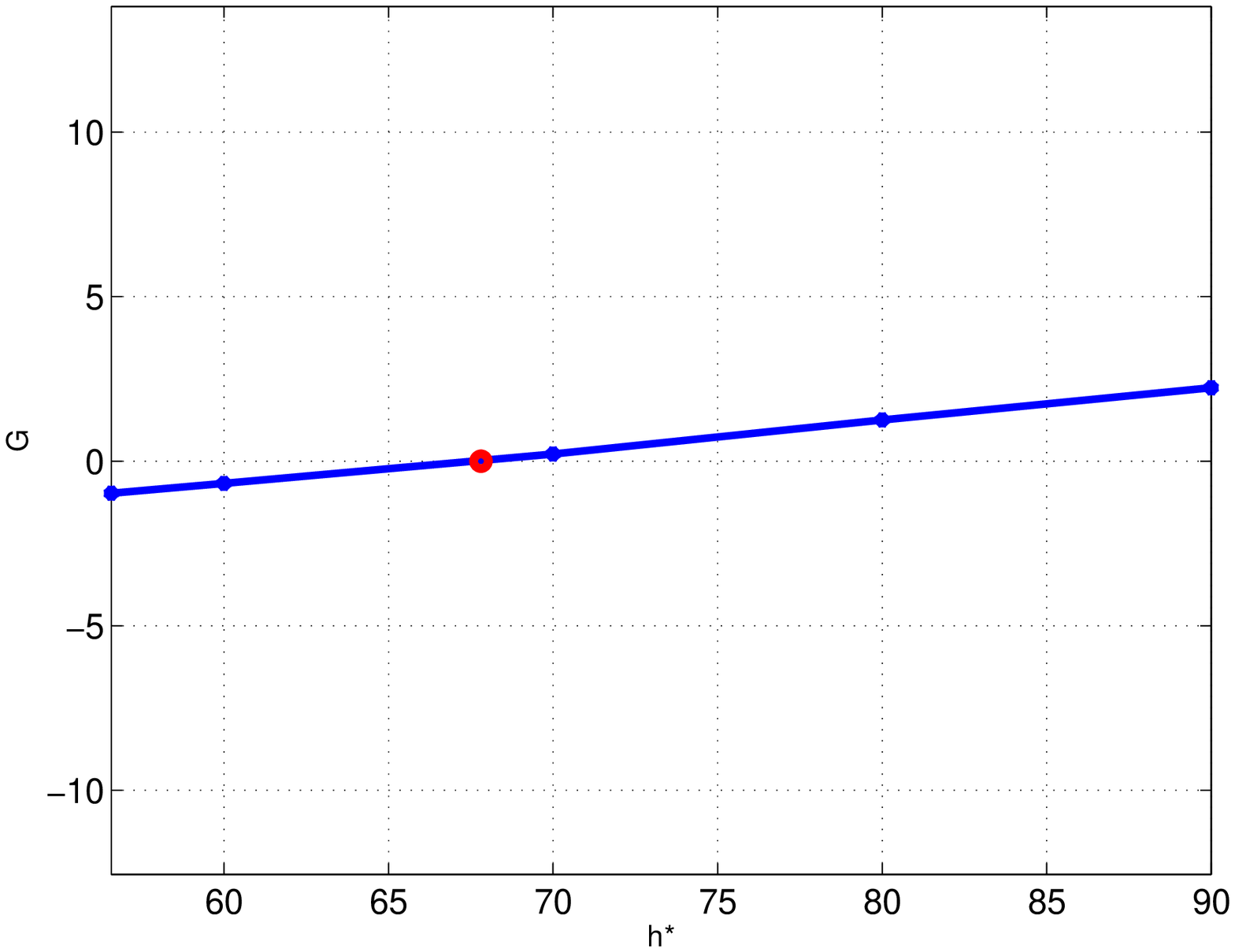}
\caption{Two cases of the $\Gamma(h^*)$ function: left and right frames are related to normal and reverse flow solutions, respectively.} 
	\label{fig:FSHGhstar}
\end{figure}

At each iteration of the ITM we have to solve the IVP 
\begin{eqnarray}\label{eq:IVP2}
& {\displaystyle \frac{d^{3}f^*}{d\eta^{*3}}} + f^* 
{\displaystyle \frac{d^{2}f^*}{d\eta^{*2}}} + \beta \left[ h_j^{*} - \left({\displaystyle
\frac{df^*}{d\eta^*}}\right)^2 \right] = 0 \nonumber \\[-1.2ex]
& \\[-1.2ex]
& f^*(0) = {\displaystyle \frac{df^*}{d\eta^*}}(0) = 0 \ , \qquad
{\displaystyle \frac{d^2f^*}{d\eta^{*2}}}(0) = \pm 1  \ . \nonumber
\end{eqnarray}
Tables~\ref{tab:Itera1} and~\ref{tab:Itera2} list the numerical iterations obtained for a sample value of $\beta$.
We notice that we solve an IVP governed by a different differential equation for each iteration because the Falkner-Skan equation is not invariant under every scaling group of point transformation.
We have chosen $ \beta=-0.01$ since, in this case, the missing initial conditions for the normal and reverse flows are not symmetric with respect to the $\beta$ axis. 
The data listed in tables~\ref{tab:Itera1} and~\ref{tab:Itera2} have been obtained by solving the modified Falkner-Skan model on $ \eta^* \in [0, 20] $ by setting 
\[
 \frac{d^2f^*}{d\eta^{*2}}(0) = \pm 1 \ ,
\]
respectively.
In both cases, we achieved convergence of the numerical results within seven iterations.
Let us now investigate the behaviour of the transformation function.
Figure~\ref{fig:FSHGhstar} shows $\Gamma(h^*)$ with respect to $h^*$ for the two cases reported in these tables.
The unique zero of the transformation function is marked by a circle.
It is worth noticing that the same scale has been used for both axes.
As it is easily seen, in both cases, we have a monotone increasing function.
We notice on the left frame, corresponding to a normal flow, that the tangent to the $\Gamma$ function at its unique zero and the $h^*$ axis define a large angle.
This is important from a numerical viewpoint because in such a case we face a well-conditioned problem.
On the other hand, this is not the case for the function plotted on the right frame of the same figure. 
The meaning is clear, reverse flow solutions are more challenging to compute than normal flow ones.
Therefore, one has to put some care when choosing the convergence criteria for the root-finder method.

\begin{figure}[p]
	\centering
\psfrag{e}[][]{\small $ \eta $} 
\psfrag{f}[][]{\small $ f(\eta)$} 
\psfrag{fe}[][]{$ {\displaystyle \frac{df}{d\eta}} $} 
\psfrag{fee}[][]{$ {\displaystyle \frac{d^2f}{d\eta^2}} $} 

\includegraphics[width=.7\textwidth]{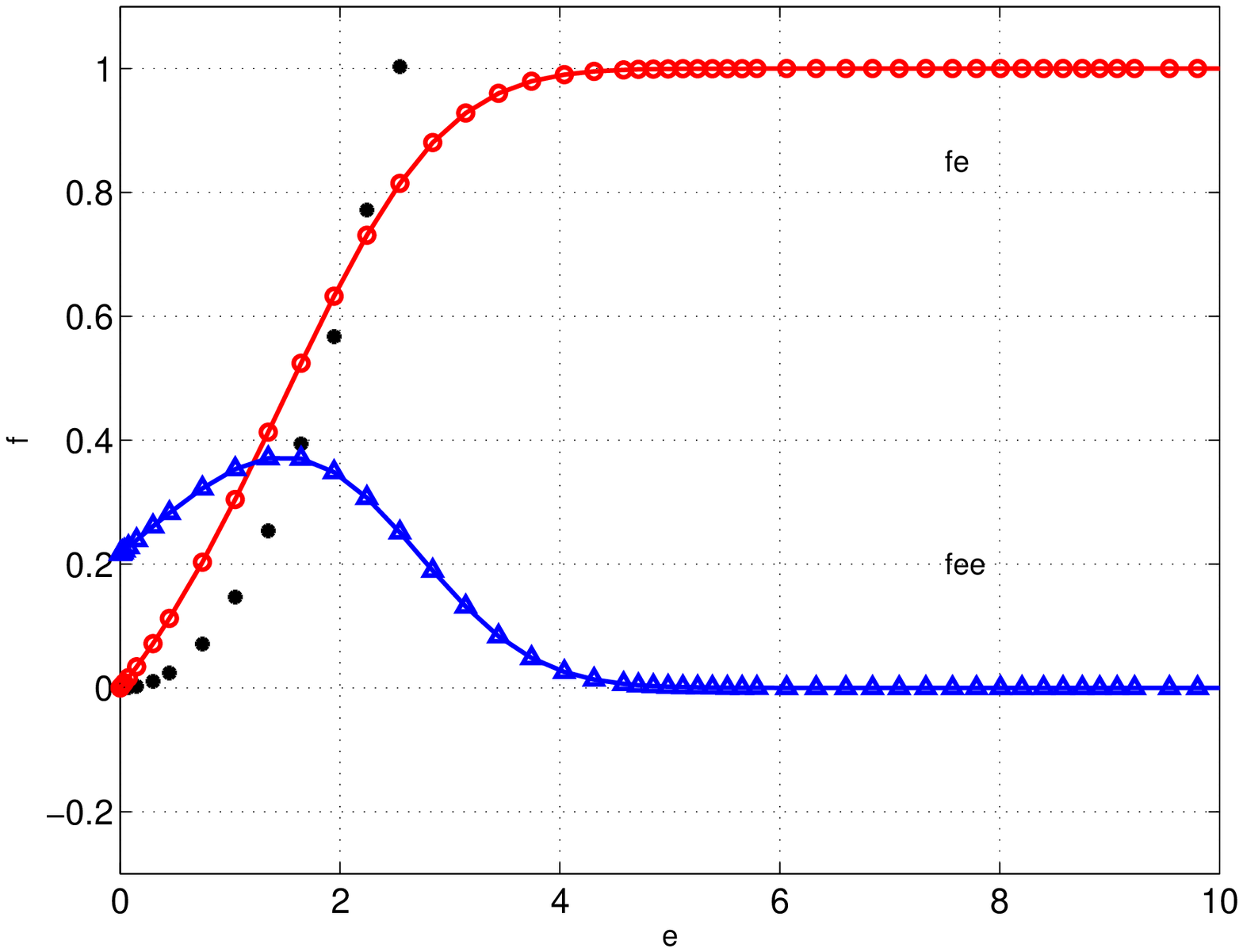} \\
\includegraphics[width=.7\textwidth]{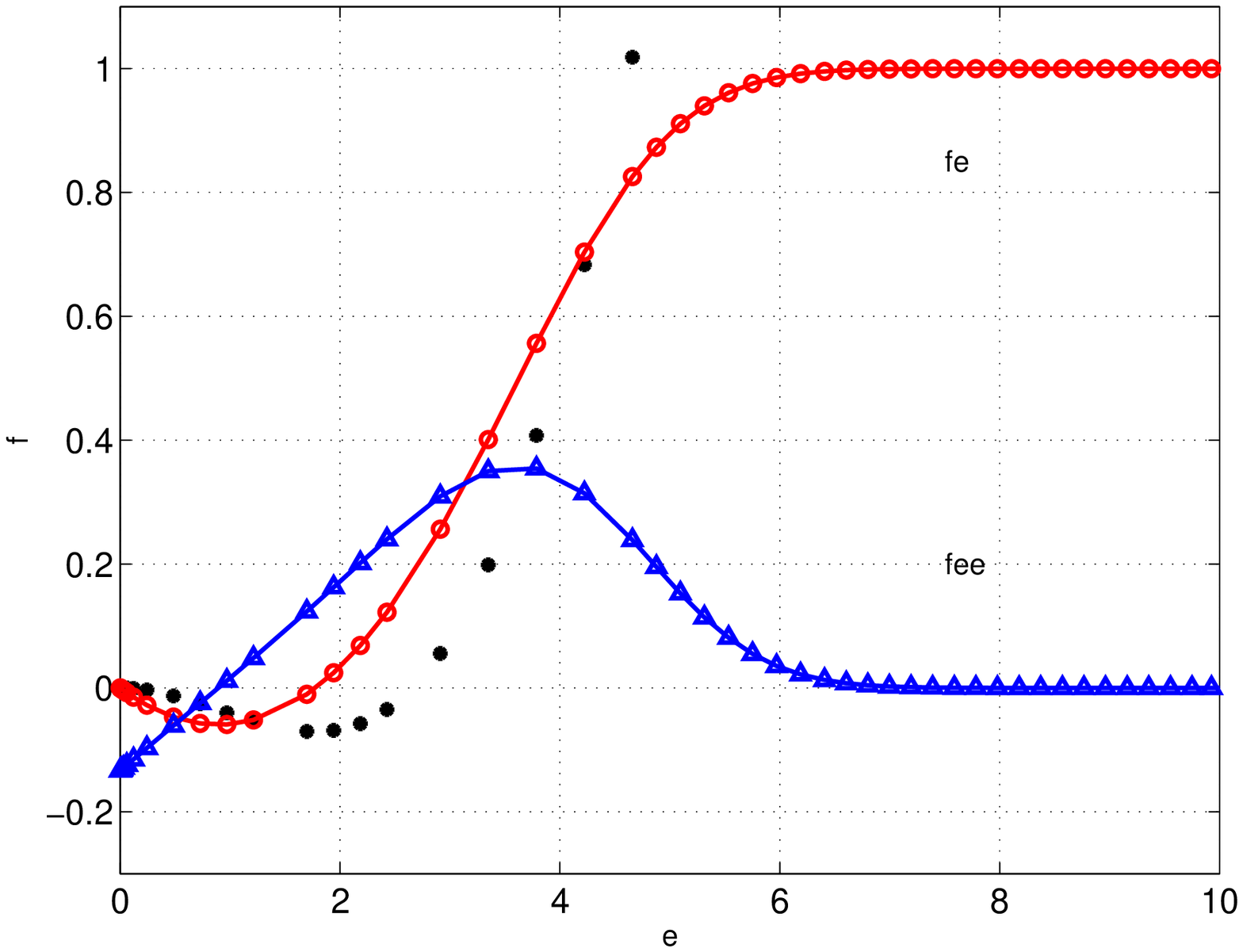}
\caption{Normal and reverse flow solutions to Falkner-Skan model for $\beta =-1.5$.
The symbols $\bullet$ denote values of $f(\eta)$.} 
	\label{fig:Falkner}
\end{figure}

Figure~\ref{fig:Falkner} shows the results of the two numerical solutions for a different value 
of $\beta$, namely $\beta = -0.15$.
In the top frame we have the normal flow and in the bottom frame we display the reverse flow solution.
In both cases the solutions were computed by introducing a truncated boundary and solving the IVP in the starred variables on $\eta^* \in [0, 20]$ with $h_0^*=1$, $h_1^* = 5$ in the top frame and
$h_0^*=15$, $h_1^* = 25$ in the bottom frame.
In this case, we achieved convergence of the numerical results within eight and seven iterations, respectively.
For the sake of clarity, we omit to plot the solutions in the starred variables computed during the iterations.
Moreover, we display only $\eta \in [0, 10]$.

In~table \ref{tab:efficiency} we report data concerning the computational cost of the IVP solver for the case $\beta = -0.15$.
Here {\it steps} is the number of successful steps, {\it failed} is the number of failed steps, and {\it evaluations} are the calculated function evaluations.

As far as the reverse flow solutions are concerned, in table~\ref{tab:Compare} we compare the missing initial condition computed by the ITM for several values of $\beta$ with results available in literature.
The agreement is really good.
It is remarkable that among the studies quoted in the introduction only a few report data related to the reverse flow solutions. 

\begin{figure}[p]
	\centering
\psfrag{b}[][]{$ \beta $} 
\psfrag{Blasius}[l][l]{Blasius} 
\psfrag{fee}[][]{${\displaystyle \frac{d^2f}{d\eta^2}} $} 

\includegraphics[width=.7\textwidth]{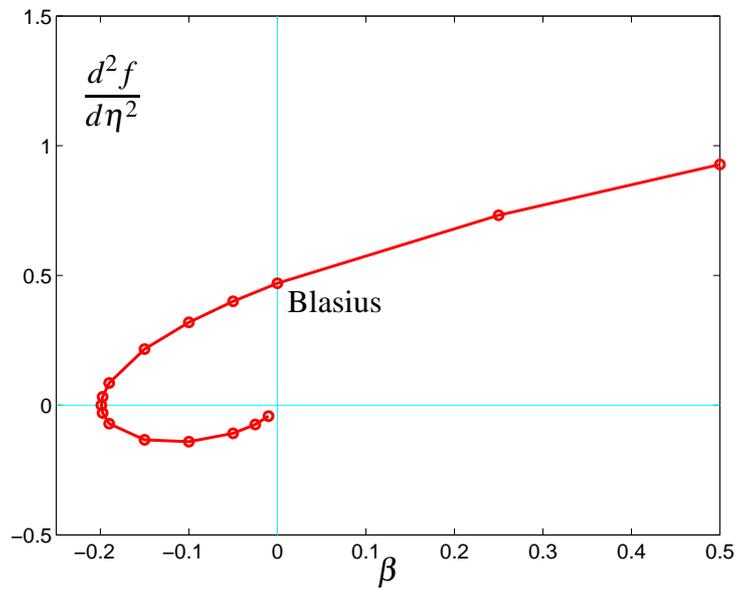}
\caption{Missing initial conditions to Falkner-Skan model for several values of $\beta$.
Positive values determine normal flow, and instead, negative values define reverse flow solutions.} 
	\label{fig:bfee}
\end{figure}

In figure~\ref{fig:bfee} we plot the behaviour of missed initial condition versus $\beta$.
The solution found by the data in table~\ref{tab:Itera2} is plotted in this figure, but not the one found in table~\ref{tab:Itera1} because this is very close to the Blasius solution.
A good initial choice of the initial iterates of $h^*$, for a given value of $\beta$, is obtained by employing values close to the one used in a successful attempt made for a close value of $\beta$.
It is interesting to note that, for values of $\beta < \beta_{\min}$ the ITM continued to iterate endlessly, whatever set of starting values for $h^*$ are selected.

Our extended algorithm has shown a kind of robustness because it is able to get convergence even when, for a chosen value of $h^*$, the IVP solver stops before arriving at the selected truncated boundary getting a wrong value of $\Gamma(h^*) = -1$.
On the other hand, the secant method gives an overflow error when this happens for two successive iterate of $h^*$.

The value of $\beta_{\min}$, corresponding to a separation point at $\eta=0$, can be found by the ITM by considering $\beta$ as a continuation parameter.
As we have seen, for $ \beta_{\min} < \beta < 0$ two solutions are available: a positive and a negative skin-friction coefficient, the missing initial condition, providing a normal and reverse flow solution.
For instance, when $\beta=-0.1988$ we get for the missing initial conditions the values
$0.005221$ and $-0.005158$, respectively.
Starting from this value of $\beta$ we can reduce it gradually and check whether the two missing initial conditions, the positive and negative values of $ {\frac{d^2f}{d\eta^2}(0)}$, converge to zero.  
Soon, we realize that we are forced the use the ITM to its natural limit.
In fact, we are trying to get a skin-friction coefficient close to zero rescaling a fixed non-zero value, plus or minus one in our case.
Anyway, when $\beta = -0.198837723795$ we found the skin-friction coefficients $6.61\mbox{D}-06$ and $-6.61\mbox{D}-06$ with 20 and 24 iterations, respectively.
Finally, we have noticed that, as far the guest for this limiting value of $\beta$ is concerned, we are allowed to reduce the chosen truncated boundary value, and for $\beta = -0.198837723795$ this truncated boundary was set equal to one, i.e. all IVP was solved on $[0, 1]$.

\begin{figure}[p]
	\centering
\psfrag{e}[][]{\small $ \eta $} 
\psfrag{f}[][]{\small $f(\eta)$} 
\psfrag{fe}[][]{$ {\displaystyle \frac{df}{d\eta}} $} 
\psfrag{fee}[][]{$ {\displaystyle \frac{d^2f}{d\eta^2}} $} 

\includegraphics[width=.7\textwidth]{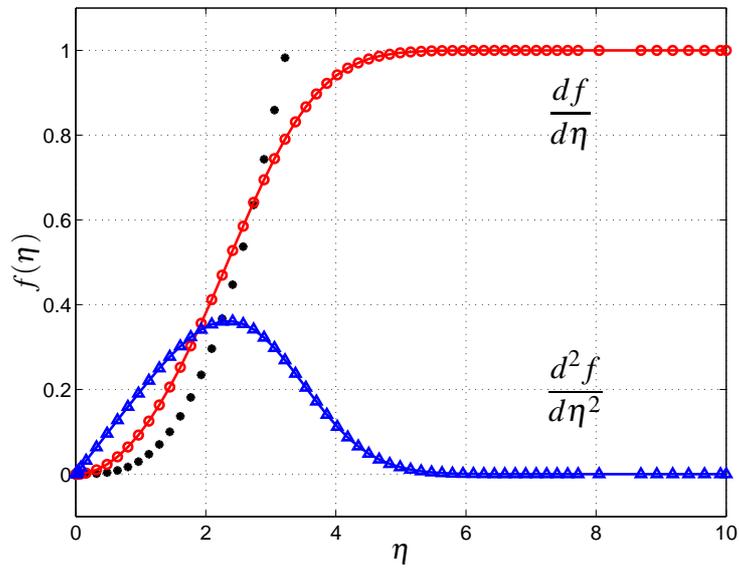}
\caption{Numerical solutions to Falkner-Skan model for $\beta = -0.1988376$. We notice that $ {\displaystyle \frac{d^2f}{d\eta^2}}(0) = 0$ and values of $f(\eta)$ are marked by $\bullet$.} 
	\label{fig:FSHL}
\end{figure}

In figure~\ref{fig:FSHL} we plot the unique solution for the limiting value $\beta_{\min}$, where $\beta_{\min}$ is given by equation (\ref{eq:bmin}).
As it easily seen this is a normal flow solution. 

The results reported so far have been found by a variable order adaptive multi-step IVP solver that was coupled up the simple secant method.
The adaptive solver uses a relative and an absolute error tolerance, for each component of the numerical solution, both equal to then to the minus six. 
As well known, the secant method is convergent provided that two initial iterates sufficiently close to the root are used, and its convergence is super-linear with an order of convergence equal to $(1+\sqrt{5})/2$.
As far as a termination criterion for the secant method is concerned, we enforced the conditions
\begin{equation}\label{eq:TC}
| \Gamma(h_j^*) | \le \mbox{Tol} \qquad \mbox{and} \qquad |h_j^*-h_{j-1}^*| \le \mbox{TolR}|h_j^*|+ \mbox{TolA} \ ,
\end{equation}
with $ \mbox{Tol}=\mbox{TolR}=\mbox{TolA}=1\mbox{D}-06$.
All computations were carried out on a 1.79 GHz AMD Turion processor with 1 GB of RAM, and for each case the execution time was few seconds.

\section{Conclusions}
In this paper we have shown how the original treatment of the Blasius problem due to T{\"o}pfer can be extended to more complex problems of boundary layer theory.
Our main concern was to solve numerically the Blasius problem, and similar problems in boundary layer theory, by initial value methods derived within scaling invariance theory.
In particular, we consider the Falkner-Skan model, for values of the parameter where multiple solutions are admitted, and report and compare the obtained numerical results, in particular data related to the famous reverse flow solutions of Stewartson.

As mentioned in the introduction, the Falkner-Skan model is a classical example where the simple shooting method cannot be applied satisfactorily for all values of the parameter $\beta$.
In fact, by applying the usual shooting method to the Falkner-Skan model, one frequently obtains
floating-point overflows in the calculations (see, for instance, Asaithambi \cite{Asaithambi:1997:NMS,Asaithambi:2005:SFE}).
A further difficulty is that the initial estimate of
the missing initial condition must occasionally be very close to the exact value in order to get convergence, cf. Cebeci and Keller \cite{Cebeci:1971:SPS}

The ITM has the same conceptual simplicity of the simple shooting method.
It is an initial value method even if we solve a different model for each iteration when the governing differential equation is not invariant under every scaling group of point transformation.
Its versatility has been shown by solving several problems of interest: 
free boundary problems \cite{Fazio:1990:SNA,Fazio:1991:ITM,Fazio:1997:NTE,Fazio:1998:SAN},
a hyperbolic moving boundary problem \cite{Fazio:1992:MBH},
the Homann and the Hiemenz flows governed by the Falkner-Skan equation in \cite{Fazio:1994:FSE},
one-dimensional parabolic moving boundary problems \cite{Fazio:2001:ITM}, two variants of the Blasius problem \cite{Fazio:2009:NTM}, namely: a boundary layer problem over moving surfaces, studied first by Klemp and Acrivos \cite{Klemp:1972:MBL}, and a boundary layer problem with slip boundary condition, that has found application to the study of gas and liquid flows at the micro-scale regime \cite{Gad-el-Hak:1999:FMM,Martin:2001:BBL}, parabolic problems on unbounded domains \cite{Fazio:2010:MBF} and, recently, see the preprints: \cite{Fazio:2012:SII} parabolic moving boundary problems, and \cite{Fazio:2012:ITM} an interesting problem in boundary layer theory: the so-called Sakiadis problem \cite{Sakiadis:1961:BLBa,Sakiadis:1961:BLBb}.

In particular, in \cite{Fazio:2001:ITM} the ITM is used to solve the sequence of free boundary problems obtained by a semi-discretization of 1D parabolic moving boundary problems.
In \cite{Fazio:2012:SII} a class of parabolic moving boundary problems is transformed to free boundary problems governed by ordinary differential equations that can be solved by the ITM.
And in \cite{Fazio:2010:MBF} a free boundary formulation for the reduced similarity models, that can be solved by the ITM, is used in order to propose a moving boundary formulation for parabolic problems on unbounded domains. 

\bigskip
\bigskip

\noindent
{\bf Acknowledgements.} The author is grateful to an anonymous reviewer for drawing his interest to the reverse flow solutions of Stewartson.
This work was supported by the University of Messina.

\begin{table}[p]
\caption{Iterations for $\beta=-0.01$ with ${\ds \frac{d^2f^*}{d\eta^{*2}}(0)} = 1$. Here and in the following the $\mbox{D}-k = 10^{-k}$ means a double precision arithmetic.}
\vspace{.5cm}
\renewcommand\arraystretch{1.3}
	\centering
		\begin{tabular}{cr@{.}lr@{.}lcr@{.}l}
\hline 
{$j$} &
\multicolumn{2}{c}%
{$h_j^*$}
& \multicolumn{2}{c}%
{$\Gamma(h_j^*)$}
&
{${\ds \frac{|h_j^*-h_{j-1}^*|}{|h_j^*|}}$}
& \multicolumn{2}{c}%
{${\displaystyle \frac{d^2f}{d\eta^2}(0)}$} \\[1.2ex]
\hline
0 &  5 &            &     0 & 631459 & & 0 & 431723 \\
1 & 10 &            &      1 & 791425  & & 0 & 384034 \\
2 &  2 & 278111 & $-$0 & 182888 & 3.389602 & 0 & 454658 \\
3 &  2 & 993420 &     0 & 0465208 & 0.238960 & 0 & 454658 \\
4 &  2 & 848366 &     9 & 5$\mbox{D}-04$ & 0.050925 & 0 & 456418 \\
5 &  2 & 845340 & $-$5 & 0$\mbox{D}-06$ & 0.001064 & 0 & 456455 \\
6 &  2 & 845356 &     6 & 1$\mbox{D}-08$ & 5.6$\mbox{D}-06$ & 0 & 456455 \\
7 &  2 & 845355 &     7 & 3$\mbox{D}-10$ & 6.7$\mbox{D}-08$ & 0 & 456455 \\
\hline			
		\end{tabular}
	\label{tab:Itera1}
\end{table}

\begin{table}[p]
\caption{Iterations for $\beta=-0.01$ with ${\ds \frac{d^2f^*}{d\eta^{*2}}(0)} = -1$.}
\vspace{.5cm}
\renewcommand\arraystretch{1.3}
	\centering
		\begin{tabular}{cr@{.}lr@{.}lcr@{.}l}
\hline 
{$j$} &
\multicolumn{2}{c}%
{$h_j^*$}
& \multicolumn{2}{c}%
{$\Gamma(h_j^*)$}
&
{${\ds \frac{|h_j^*-h_{j-1}^*|}{|h_j^*|}}$}
& \multicolumn{2}{c}%
{${\displaystyle \frac{d^2f}{d\eta^2}(0)}$} \\[1.2ex]
\hline
0 &  75 &            &     0 & 731890 & & $-$0 & 059237 \\
1 & 150 &           &     5 & 263092 & & $-$0 & 092368 \\
2 &  62 & 885833 & $-$0 & 443040 & 1.385275 & $-$0 & 028870 \\
3 &  69 & 649620 &     0 & 181067 & 0.097112 & $-$0 & 046991 \\
4 & 67 & 687299 & $-$0 & 011297  & 0.028991 & $-$0 & 042016 \\
5 & 67 & 802542 & $-$2 & 1$\mbox{D}-04$ & 0.001700 & $-$0 & 042315 \\
6 & 67 & 804749 & 2 & 8$\mbox{D}-07$      & 3.3$\mbox{D}-05$ & $-$0 & 042321 \\
7 & 67 & 804746 & 7 & 9$\mbox{D}-10$    & 4.3$\mbox{D}-08$ & $-$0 & 042321 \\
\hline			
		\end{tabular}
	\label{tab:Itera2}
\end{table}

\begin{table}[p]
\caption{Efficiency of the ITM for $\beta = -0.15$ and reverse flow iterations.}
\vspace{.5cm}
\renewcommand\arraystretch{1.3}
	\centering
		\begin{tabular}{cccc}
\hline 
{$j$} & {\it steps} & {\it failed} & {\it evaluations} \\
\hline
$0$ & 447 & 47 & ~942\\
$1$ & 566 & 84 & 1217\\
$2$ & 458 & 47 & ~964\\
$3$ & 483 & 47 & 1014\\
$4$ & 464 & 52 & ~981\\
$5$ & 425 & 58 & ~909\\
$6$ & 479 & 49 & 1008\\
$7$ & 463 & 46 & ~973\\
\hline			
		\end{tabular}
	\label{tab:efficiency}
\end{table}

\begin{table}[p]
\caption{Comparison for the reverse flow skin-friction coefficients ${\displaystyle \frac{d^2f}{d\eta^{2}}(0)}$.
For all cases we used $h_0^* = 15$ and $h_1^* = 25$.
The iterations were, from top to bottom line: $8$, $7$, $9$, $7$, and $7$.}
\vspace{.5cm}
\renewcommand\arraystretch{1.3}
	\centering
		\begin{tabular}{r@{.}lcccc}
\hline 
\multicolumn{2}{c}%
{$\beta$} &
Stewartson \cite{Stewartson:1954:FSF} &
Asaithambi \cite{Asaithambi:1997:NMS} & Auteri et al. \cite{Auteri:2012:GLS}&
ITM \\[1ex]
\hline
$-0$ & $025$ & $-0.074$ & & $$ & $-0.074366$ \\
$-0$ & $05$  & $-0.108$ & & $$ & $-0.108271$ \\
$-0$ & $1$    & $-0.141$ & $-0.140546$ & $-0.140546$ & $-0.140546$ \\
$-0$ & $15$  & $-0.132$ & $-0.133421$ & $-0.133421$ & $-0.133421$ \\
$-0$ & $18$  & $-0.097$ & $-0.097692$ & $-0.097692$ & $-0.097692$ \\
\hline			
		\end{tabular}
	\label{tab:Compare}
\end{table}

\end{document}